\documentclass[12pt,fleqn,leqno]{article}

\usepackage{amsmath,amssymb}
\usepackage{float}
\usepackage[dvips]{graphicx}
\usepackage[dvips]{color}

\numberwithin{equation}{section}

\makeatletter
\renewcommand{\section}{\@startsection{section}{1}{0pt}{20pt}{6pt}{\large\bf}}
\renewcommand{\@seccntformat}[1]{\csname the#1\endcsname.\ }

\def\footnoterule{\kern -3pt \hrule width 2.7 true cm \kern 2.6pt}

\def\h{\hspace}

\def\ni{\noindent}
\def\p{\!+\!}
\def\m{\!-\!}

\def\EE{\mathsf E\:\!}
\def\PP{\mathsf P}

\def\cF{{\cal F}}

\def\R{I\!\!R}
\def\LL{I\!\!L}
\def\eps{\varepsilon}

\begin{document}

\title{\bf Quickest Real-Time Detection of \\ a Brownian Coordinate Drift}
\author{P. A. Ernst \& G. Peskir}
\date{}
\maketitle




{\par \leftskip=2cm \rightskip=2cm \footnotesize

Consider the motion of a Brownian particle in two or more
dimensions, whose coordinate processes are standard Brownian motions
with zero drift initially, and then at some random/unobservable
time, one of the coordinate processes gets a (known) non-zero drift
permanently. Given that the position of the Brownian particle is
being observed in real time, the problem is to detect the time at
which a coordinate process gets the drift as accurately as possible.
We solve this problem in the most uncertain scenario when the
random/unobservable time is (i) exponentially distributed and (ii)
independent from the initial motion without drift. The solution is
expressed in terms of a stopping time that minimises the probability
of a false early detection and the expected delay of a missed late
detection. To our knowledge this is the first time that such a
problem has been solved exactly in the literature.

\par}


\footnote{{\it Mathematics Subject Classification 2010.} Primary
60G40, 60J65, 60H30. Secondary 35J15, 45G10, 62C10.}

\footnote{{\it Key words and phrases:} Quickest detection, Brownian
motion, optimal stopping, elliptic partial differential equation,
free-boundary problem, smooth fit, nonlinear Fredholm integral
equation, the change-of-variable formula with local time on
surfaces.}

\vspace{-6mm}

\section{Introduction}

Imagine the motion of a Brownian particle in two or more dimensions,
whose coordinate processes are standard Brownian motions with zero
drift initially, and then at some random/ unobservable time
$\theta$, one of the coordinate processes gets a (known) non-zero
drift $\mu$ permanently. Assuming that the position of the Brownian
particle is being observed in real time, the problem is to detect
the time $\theta$ at which a coordinate process gets the drift $\mu$
as accurately as possible. The purpose of the present paper is
derive the solution to this problem in the most uncertain scenario
when $\theta$ is assumed to be (i) exponentially distributed and
(ii) independent from the initial motion without drift.

Denoting the position of the Brownian particle in two or more
dimensions by $X$, the error to be minimised over all stopping times
$\tau$ of $X$ is expressed as the the linear combination of the
probability of the \emph{false alarm} $\PP_{\!\pi}(\tau\! <\!
\theta)$ and the expected \emph{detection delay} $\EE_\pi(\tau \m
\theta)^+$ where $\pi \in [0,1]$ denotes the probability that
$\theta$ has already occurred at time $0$. This problem formulation
of quickest detection dates back to \cite{Sh-1} and has been
extensively studied to date (see \cite{Sh-3} and the references
therein). The linear combination represents the Lagrangian and once
the optimal stopping problem has been solved in this form it will
also lead to the solution of the constrained problems where an upper
bound is imposed on either the probability of the false alarm or the
expected detection delay respectively.

A canonical example is the standard Brownian motion in one dimension
with one \emph{constant} drift changing to another. This problem has
also been solved in finite horizon (see \cite{GP} and the references
therein). Books \cite[Section 4.4]{Sh-2} and \cite[Section 22]{PS}
contain expositions of these results and provide further details and
references. The \emph{signal-to-noise ratio} (defined as the
difference between the new drift and the old drift divided by the
diffusion coefficient) in all these problems is \emph{constant} so
that the resulting optimal stopping problem for the posterior
probability distribution ratio process $\varPhi$ of $\theta$ given
$X$ is \emph{one-dimensional}. A more general problem formulation
for diffusion processes $X$ in one dimension when one
\emph{non-constant} drift changes to another has been considered in
\cite{GS}. A specific problem of this kind when $X$ is a Bessel
process has been solved in \cite{JP}. The signal-to-noise ratio in
these problems is not constant and the resulting optimal stopping
problem for $\varPhi$ coupled with $X$ (to make it Markovian) is
\emph{two-dimensional}. The infinitesimal generator of the
Markov/diffusion process $(\varPhi,X)$ in these problems is of
\emph{parabolic} type.

Related quickest detection problems for $X$ in two dimensions have
been studied in \cite{BP} and \cite{DPS}. The change of
probabilistic characteristics in these problems can affect both
coordinate processes of $X$ and not only one as in the present
paper. The coordinate processes of $X$ in \cite{BP} are Poisson
processes and the resulting two-dimensional optimal stopping problem
for $\varPhi$ has been studied using an iteration technique. The
coordinate processes in \cite{DPS} are Wiener/Poisson processes and
the resulting optimal stopping problem for $\varPhi$ is
one-dimensional.

The quickest detection setting of the observed process $X$ in two or
more dimensions may also be viewed as a \emph{multi-channel sensor
system}. Quickest detection problems of this kind in two or more
dimensions have been studied in a number of papers (see \cite{ZRH}
\& \cite{FS} and the references therein). These papers usually
establish `asymptotic optimality' of an `ad-hoc' stopping rule and
no `exact' (optimal) solution has been derived in the literature to
date.

In contrast to the quickest detection problems solved to date, we
will see below that the \emph{multi-dimensional} Markov/diffusion
process $\varPhi$ in the quickest detection problem of the present
paper has the infinitesimal generator of \emph{elliptic} type.
Finding the \emph{exact} solution to the quickest detection problem
for the observed process $X$ in \emph{two or more} dimensions is the
main contribution of the present paper. To our knowledge this is the
first time that such a problem has been solved exactly in the
literature.

\section{Formulation of the problem}

In this section we formulate the quickest detection problem under
consideration. The initial formulation of the problem will be
revaluated under a change of measure in the next section. To
simplify the exposition we will assume throughout that the observed
process is \emph{two-dimensional}. This assumption will be extended
to \emph{three or more dimensions} in the final section below.

\vspace{4mm}

1.\ We consider a Bayesian formulation of the problem where it is
assumed that one observes a sample path of the standard
two-dimensional Brownian motion $X=(X^1,X^2)$, whose coordinate
processes $X^1$ and $X^2$ are standard Brownian motions with zero
drift initially, and then at some random/unobservable time $\theta$
taking value $0$ with probability $\pi \in [0,1]$ and being
exponentially distributed with parameter $\lambda>0$ given that
$\theta>0$, one of the coordinate processes $X^1$ and $X^2$ gets a
(known) non-zero drift $\mu$ permanently. The problem is to detect
the time $\theta$ at which a coordinate process gets the drift $\mu$
as accurately as possible (neither too early nor too late). This
problem belongs to the class of quickest real-time detection
problems as discussed in Section 1 above.

\vspace{4mm}

2.\ The observed process $X=(X^1,X^2)$ solves the stochastic
differential equations
\begin{align} \h{7pc} \label{2.1}
&dX_t^1 = \mu\:\! I(\beta\! =\! 1,t\! \ge\! \theta)\;\! dt +
dB_t^1 \\[2pt] \label{2.2} &dX_t^2 = \mu\:\! I(\beta\! =\! 2,t
\! \ge\! \theta) \;\! dt + dB_t^2
\end{align}
driven by a standard two-dimensional Brownian motion $B=(B^1,B^2)$
under the probability measure $\PP_{\!\pi}$ specified below, where
the random variable $\beta$ satisfies $\PP_{\!\pi}(\beta\! =\!
1)=p_1$ and $\PP_{\!\pi}(\beta\! =$ $\! 2)=p_2$ for some $p_1,p_2
\in [0,1]$ with $p_1 \p p_2 = 1$ given and fixed, meaning that
$\beta=i$ if and only if the coordinate process $X_i$ gets drift
$\mu$ at time $\theta$ with probability $p_i$ for $i=1,2$. The
unobservable time $\theta$, the unknown coordinate $\beta$, and the
driving Brownian motion $B$ are all assumed to be independent under
$\PP_{\!\pi}$ for $\pi \in [0,1]$ given and fixed.

\vspace{4mm}

3.\ Standard arguments imply that the previous setting can be
realised on a probability space $(\Omega,\cF,\PP_{\!\pi})$ with the
probability measure $\PP_{\!\pi}$ being decomposable as follows
\begin{equation} \h{1pc} \label{2.3}
\PP_{\!\pi} = p_1\;\! \pi\;\! \PP_{\!1}^0 +  p_2\;\! \pi\;\!
\PP_{\!2}^0 + p_1\;\! (1 \m \pi)\! \int_0^\infty\! \lambda\;
\! e^{-\lambda t}\, \PP_{\!1}^t\, dt + p_2\;\! (1 \m \pi)\!
\int_0^\infty\! \lambda\;\! e^{-\lambda t}\, \PP_{\!2}^t\, dt
\end{equation}
for $\pi \in [0,1]$ where $\PP_{\!i}^t$ is the probability measure
under which the coordinate process $X^i$ gets drift $\mu$ at time $t
\in [0,\infty)$ for $i=1,2$. The decomposition \eqref{2.3} expresses
the fact that the unobservable time $\theta$ is a non-negative
random variable satisfying $\PP_{\!\pi}(\theta = 0) = \pi$ and
$\PP_{\!\pi}(\theta > t\, \vert\, \theta > 0) = e^{-\lambda t}$ for
$t>0$. Thus $\PP_{\!i}^t(X \in\, \cdot\, )$ $= \PP_{\!\pi}(X \in\,
\cdot\; \vert\, \beta=i,\:\! \theta=t)$ is the probability law of
the standard two-dimensional Brownian motion process $X=(X^1,X^2)$
whose coordinate process $X^i$ gets drift $\mu$ at time $t \in
[0,\infty)$ for $i=1,2$. To remain consistent with this notation we
also denote by $\PP_{\!i}^\infty$ the probability measure under
which the coordinate process $X^i$ of the observed process
$X=(X^1,X^2)$ gets no drift $\mu$ at a finite time for $i=1,2$. Thus
$\PP_{\!i}^\infty(X \in\, \cdot\; ) = \PP_{\!\pi}(X \in\, \cdot\;
\vert\, \beta=i,\:\! \theta = \infty)$ is the probability law of the
standard two-dimensional Brownian motion process for $i=1,2$.
Clearly the subscript $i$ is superfluous in this case and we will
often write $\PP^\infty$ instead of $\PP_{\!i}^\infty$ for $i=1,2$.
Moreover, by $\PP_{\!i}$ we denote the probability measure under
which the coordinate process $X^i$ gets drift $\mu$ at time $\theta$
for $i=1,2$. From \eqref{2.3} we see that
\begin{equation} \h{9pc} \label{2.4}
\PP_{\!\pi} = p_1\:\! \PP_{\!1} +  p_2\:\! \PP_{\!2}
\end{equation}
where $\PP_{\!i} = \pi\;\! \PP_{\!i}^0 + (1 \m \pi)\!
\int_0^\infty\! \lambda\; \! e^{-\lambda t}\, \PP_{\!i}^t\, dt$ for
$i=1,2$ and $\pi \in [0,1]$. Note that $\PP_{\!i}$ depends on $\pi
\in [0,1]$ as well but we will omit this dependence from its
notation for $i=1,2$.

\vspace{2mm}

4.\ Being based upon continuous observation of $X=(X^1,X^2)$, the
problem is to find a stopping time $\tau_*$ of $X$ (i.e.\ a stopping
time with respect to the natural filtration $\cF_t^X = \sigma(X_s\,
\vert\, 0 \le s \le t)$ of $X$ for $t \ge 0$) that is `as close as
possible' to the unknown time $\theta$. More precisely, the problem
consists of computing the value function
\begin{equation} \h{5.5pc} \label{2.5}
V(\pi) = \inf_\tau \Big[ \PP_{\!\pi}(\tau < \theta) + c\;\! \EE_\pi
(\tau - \theta)^+ \Big]
\end{equation}
and finding the optimal stopping time $\tau_*$ at which the infimum
in \eqref{2.5} is attained for $\pi \in [0,1]$ and $c>0$ given and
fixed (recalling also that $p_1,p_2 \in [0,1]$ with $p_1 \p p_2 = 1$
are given and fixed). Note in \eqref{2.5} that $\PP_{\!\pi}(\tau <
\theta)$ is the probability of the \emph{false alarm} and
$\EE_\pi(\tau - \theta)^+$ is the expected \emph{detection delay}
associated with a stopping time $\tau$ of $X$ for $\pi \in [0,1]$.
Recall also that the expression on the right-hand side of
\eqref{2.5} is the Lagrangian associated with the constrained
problems as discussed in Section 1 above.

\vspace{4mm}

5.\ To tackle the optimal stopping problem \eqref{2.5} we consider
the \emph{posterior probability distribution process}
$\varPi=(\varPi_t)_{t \ge 0}$ of $\theta$ given $X$ that is defined
by
\begin{equation} \h{8pc} \label{2.6}
\varPi_t = \PP_{\!\pi}(\theta \le t\, \vert\, \cF_t^X)
\end{equation}
for $t \ge 0$. Note that we have
\begin{equation} \h{8pc} \label{2.7}
\varPi_t = \varPi_t^1 + \varPi_t^2
\end{equation}
where we set
\begin{equation} \h{2pc} \label{2.8}
\varPi_t^1 = \PP_{\!\pi}(\beta=1,\: \theta \le t\, \vert\, \cF_t^X)
\;\;\; \&\;\;\; \varPi_t^2 = \PP_{\!\pi}(\beta=2,\: \theta \le t\,
\vert \, \cF_t^X)
\end{equation}
for $t \ge 0$. The right-hand side of \eqref{2.5} can be rewritten
to read
\begin{equation} \h{5pc} \label{2.9}
V(\pi) = \inf_\tau\;\! \EE_\pi \Big( 1 \m \varPi_\tau + c \int_0^\tau\!
\varPi_t\, dt \Big)
\end{equation}
for $\pi \in [0,1]$.

\vspace{6mm}

6.\ To connect the process $\varPi$ to the observed process $X$ we
set
\begin{equation} \h{2pc} \label{2.10}
\bar \varPi_t^1 = \PP_{\!\pi}(\beta=1,\: \theta > t\, \vert\, \cF_t^X)
\;\;\; \&\;\;\; \bar \varPi_t^2 = \PP_{\!\pi}(\beta=2,\: \theta > t\,
\vert \, \cF_t^X)
\end{equation}
and define the \emph{posterior probability distribution ratio
process} $\varPhi=(\varPhi^1,\varPhi^2)$ of $\theta$ given $X$ by
\begin{equation} \h{7pc} \label{2.11}
\varPhi_t^1 = \frac{\varPi_t^1}{\bar \varPi_t^1}\;\;\; \&\;\;\;
\varPhi_t^2 = \frac{\varPi_t^2}{\bar \varPi_t^2}
\end{equation}
for $t \ge 0$. Using \eqref{2.3} we find that
\begin{equation} \h{5pc} \label{2.12}
\varPi_t^i = p_i\;\! \pi\;\! \frac{d \PP_{\!i,t}^0}{d \PP_{\!\pi,
t}} + p_i\;\! (1 \m \pi)\! \int_0^t\! \lambda\;\! e^{-\lambda s}
\, \frac{d \PP_{\!i,t}^s}{d \PP_{\!\pi,t}}\, ds
\end{equation}
where $\PP_{\!i,t}^s$ and $\PP_{\!\pi,t}$ denote the restrictions of
the measures $\PP_{\!i}^s$ and $\PP_{\!\pi}$ to $\cF_t^X$ for $s \ge
0$ and $i=1,2$ respectively. Similarly, using \eqref{2.3} we find
that
\begin{equation} \h{7pc} \label{2.13}
\bar \varPi_t^i = p_i\;\! (1 \m \pi)\;\! e^{-\lambda t}\;\! \frac
{d \PP_{\!t}^\infty}{d \PP_{\!\pi,t}}
\end{equation}
where $\PP_{\!t}^\infty$ and $\PP_{\!\pi,t}$ denote the restrictions
of the measures $\PP^\infty$ and $\PP_{\!\pi}$ to $\cF_t^X$ for $t
\ge 0$ and $i=1,2$ (notice in this derivation that $d
\PP_{\!i,t}^s/d \PP_{\!\pi,t} = d \PP_{t}^\infty/d \PP_{\!\pi,t}$
for $s \ge t$). From \eqref{2.12} and \eqref{2.13} we see that
taking ratios as in \eqref{2.11} removes dependence on
$\PP_{\!\pi,t}$ which makes explicit calculations possible. Indeed,
using the Girsanov theorem we see that the \emph{likelihood ratio
process} $L=(L^1,L^2)$ can be expressed as follows
\begin{equation} \h{7pc} \label{2.14}
L_t^i = \frac{d \PP_{\!i,t}^0}{d \PP_{\!t}^\infty} = \exp \Big(
\mu\:\! X_t^i - \frac{\mu^2}{2}\:\! t \Big)
\end{equation}
for $t \ge 0$ and $i=1,2$. Moreover, using \eqref{2.12} and
\eqref{2.13} we find by \eqref{2.11} that
\begin{equation} \h{7pc} \label{2.15}
\varPhi_t^i = e^{\lambda t} L_t^i\:\! \Big( \varPhi_0^i + \lambda
\int_0^t \frac{ds} {e^{\lambda s} L_s^i}\, \Big)
\end{equation}
with $\varPhi_0^i=\pi/(1 \m \pi)$ for $t \ge 0$ and $i=1,2$
(\:\!notice in this derivation that $d \PP_{\!i,t}^s/d \PP_{\!i,t}^t
= L_t^i/L_s^i$ for $s \le t$). From \eqref{2.14} and \eqref{2.15} we
see that the process $\varPhi=(\varPhi^1,\varPhi^2)$ is an explicit
(path-dependent) functional of the observed process $X=(X^1,X^2)$
and hence observable (by observing a sample path of $X$ we are also
seeing a sample path of $\varPhi$ both in real time).

\section{Measure change}

In this section we show that changing the probability measure
$\PP_{\!\pi}$ for $\pi \in [0,1]$ to $\PP^\infty$ in the optimal
stopping problem \eqref{2.5} or \eqref{2.9} provides crucial
simplifications of the setting which make the subsequent analysis
possible. This will be achieved by invoking the decomposition of
$\PP_{\!\pi}$ into $\PP_{\!1}$ and $\PP_{\!2}$ as stated in
\eqref{2.4} above, changing both probability measures $\PP_{\!1}$
and $\PP_{\!2}$ to $\PP_{\!1}^\infty$ and $\PP_{\!2}^\infty$
respectively, and recalling that both $\PP_{\!1}^\infty$ and
$\PP_{\!2}^\infty$ coincide with $\PP^\infty$.

\vspace{4mm}

1.\ We show that the optimal stopping problem \eqref{2.9} admits a
transparent reformulation under the probability measure $\PP^\infty$
in terms of the process $\varPhi=(\varPhi^1,\varPhi^2)$ defined by
\eqref{2.11} above. Recall that $\varPhi^i$ starts at $\pi/(1 \m
\pi)$ and this dependence on the initial point will be indicated by
a superscript to $\varPhi^i$ when needed for $i=1,2$.

\vspace{6mm}

\textbf{Proposition 1.} \emph{The value function $V$ from
\eqref{2.9} satisfies the identity
\begin{equation} \h{7pc} \label{3.1}
V(\pi) = (1 \m \pi)\, \big[ 1 + c\;\! \hat V(\pi) \big]
\end{equation}
where the value function $\hat V$ is given by
\begin{equation} \h{2pc} \label{3.2}
\hat V(\pi) = \inf_\tau\;\! \EE^\infty \Big[ \int_0^\tau e^{-\lambda t}
\Big( p_1\:\! \varPhi_t^{1,\pi/(1-\pi)}\! + p_2\:\! \varPhi_t^{2,\pi/
(1-\pi)} - \frac{\lambda} {c}\;\! \Big)\;\! dt\:\! \Big]
\end{equation}
for $\pi \in [0,1)$ and the infimum in \eqref{3.2} is taken over all
stopping times $\tau$ of $X$.}

\vspace{6mm}

\textbf{Proof.} Let a (bounded) stopping time $\tau$ of $X$ be given
and fixed. Set
\begin{equation} \h{4pc} \label{3.3}
\tilde \varPi_t^1 = \PP_{\!1}(\theta \le t\, \vert\, \cF_t^X)
\;\;\; \&\;\;\; \tilde \varPi_t^2 = \PP_{\!2}(\theta \le t\,
\vert \, \cF_t^X)
\end{equation}
for $t \ge 0$. We claim that
\begin{equation} \h{-0.2pc} \label{3.4}
\EE_\pi \Big( 1 \m \varPi_\tau + c\! \int_0^\tau\! \varPi_t\,
dt \Big) = p_1\, \EE_1 \Big( 1 \m \tilde \varPi_\tau^1 + c\!
\int_0^\tau\! \tilde \varPi_t^1\, dt \Big) + p_2\, \EE_2 \Big(
1 \m \tilde \varPi_\tau^2 + c\! \int_0^\tau\! \tilde \varPi_t^2
\, dt \Big)
\end{equation}
for $\pi \in [0,1]$. For this, first note that \eqref{2.7} yields
\begin{equation} \h{-0.2pc} \label{3.5}
\EE_\pi \Big( 1 \m \varPi_\tau + c\! \int_0^\tau\! \varPi_t\,
dt \Big) = 1 + \EE_\pi \Big( \m \varPi_\tau^1 \m \varPi_\tau^2
+ c\! \int_0^\tau \! \varPi_t^1\, dt + c\! \int_0^\tau\!
\varPi_t^2\, dt \Big)
\end{equation}
for $\pi \in [0,1]$. Next note that \eqref{2.4} implies that
\begin{equation} \h{2pc} \label{3.6}
\EE_\pi \big( \Pi_\tau^i \big) = \PP_{\!\pi}(\beta=i,\: \theta
\le \tau) = p_i\;\! \PP_{\!i}(\theta \le \tau) = p_i\;\! \EE_i
\big( \tilde \Pi_\tau^i \big)
\end{equation}
and similarly we find that
\begin{align} \h{-0.2pc} \label{3.7}
\EE_\pi \Big( \int_0^\tau\! \varPi_t^i\, dt\;\! \Big) &= \EE_\pi
\Big( \int_0^\infty\! \PP_{\!\pi}(\beta=i,\: \theta \le t \le
\tau\, \vert\, \cF_t^X)\, dt\;\! \Big) = \int_0^\infty\! \PP_{
\!\pi}(\beta=i,\: \theta \le t \le \tau)\, dt \\ \notag &= p_i
\! \int_0^\infty\! \PP_{\!i}(\theta \le t \le \tau)\, dt = p_i
\, \EE_i \Big( \int_0^\infty\! \PP_{\!i}(\theta \le t \le \tau
\, \vert\, \cF_t^X)\, dt\;\! \Big) \\ \notag &= p_i\;\! \EE_i
\Big( \int_0^\tau\! \tilde \varPi_t^i\, dt\;\! \Big)
\end{align}
for $\pi \in [0,1]$ and $i=1,2$. Finally, combining
\eqref{3.5}-\eqref{3.7} we obtain \eqref{3.4} as claimed.

Focusing on each of the two expectations on the right-hand side of
\eqref{2.4} separately, and noticing that the enlargement of the
filtration from $\cF_t^{X^i}$ to $\cF_t^X$ for $t \ge 0$ creates no
difficulty for $i=1,2$ because $X^1$ and $X^2$ are independent under
both $\PP_{\!1}$ and $\PP_{\!2}$, we see that the problem of
establishing \eqref{3.1} and \eqref{3.2} reduces to one dimension.
Hence applying the change-of-measure identity (4.12) from \cite{JP}
to each of the two expectations on the right-hand side of
\eqref{2.4} separately, we obtain
\begin{equation} \h{2pc} \label{3.8}
\EE_i \Big( 1 \m \tilde \varPi_\tau^i + c\! \int_0^\tau\! \tilde
\varPi_t^i\, dt \Big) = (1 \m \pi) \Big( 1 + c\, \EE_i^\infty \Big[
\int_0^\tau e^{-\lambda t} \Big( \tilde \varPhi_t^i - \frac{\lambda}
{c}\;\! \Big)\;\! dt\:\! \Big] \Big)
\end{equation}
for $\pi \in [0,1]$ and $i=1,2$. On closer look we see that $\tilde
\varPhi^i$ and $\varPhi^i$ coincide for $i=1,2$ (which is not
surprising in view of \eqref{2.14} above). Recalling that
$\PP_{\!1}^\infty$ and $\PP_{\!2}^\infty$ coincide with $\PP^\infty$
and inserting \eqref{3.8} into \eqref{3.4} we see that \eqref{3.1}
and \eqref{3.2} hold as claimed. \hfill $\square$

\vspace{4mm}

2.\ From Proposition 1 we see that the optimal stopping problem
\eqref{2.5} or \eqref{2.9} is equivalent to the optimal stopping
problem \eqref{3.2}. Using the fact pointed out in the proof above
that $\varPhi^i$ and $\tilde \varPhi^i$ coincide for $i=1,2$, we see
from (4.7) in \cite{JP} that $\varPhi^1$ and $\varPhi^2$ solve the
following stochastic differential equations
\begin{align} \h{7pc} \label{3.9}
d \varPhi_t^1 = \lambda\:\! (1 \p \varPhi_t^1)\;\! dt + \mu\;\!
\varPhi_t^1\;\! dB_t^1 \\[3pt] \label{3.10} d \varPhi_t^2 = \lambda
\:\! (1 \p \varPhi_t^2)\;\! dt + \mu\;\! \varPhi_t^2\;\! dB_t^2
\end{align}
under $\PP^\infty$ with $\varPhi_0^1 = \varphi_1$ and $\varPhi_0^2 =
\varphi_2$ in $[0,\infty)$ both being equal to $\pi/(1 \m \pi)$ for
$\pi \in [0,1)$. The system of stochastic differential equations
\eqref{3.9}-\eqref{3.10} has a unique strong solution given by
\eqref{2.14}+\eqref{2.15} above. Hence the process $\varPhi =
(\varPhi^1,\varPhi^2)$ is both strong Markov and strong Feller (see
e.g.\ \cite[pp 158-163 \& pp 170-173]{RW}). Basic properties of the
one-dimensional diffusion processes $\varPhi^1$ and $\varPhi^2$ are
reviewed in \cite[Section 2]{Pe-2}. In particular, it is known that
$\varPhi^i$ is \emph{recurrent} in $[0,\infty)$ if and only if
$\lambda \le \mu^2\!/2$ for $i=1,2$. If $\lambda > \mu^2\!/2$ then
$\varPhi^i$ is \emph{transient} in $[0,\infty)$ with $\varPhi_t^i
\rightarrow \infty$ almost surely under $\PP^\infty$ as $t
\rightarrow \infty$ for $i=1,2$.

\vspace{6mm}

3.\ To tackle the equivalent optimal stopping problem \eqref{3.2}
for the strong Markov process $\varPhi = (\varPhi^1,\varPhi^2)$
solving \eqref{3.9}-\eqref{3.10} we will enable $\varPhi =
(\varPhi^1,\varPhi^2)$ to start at any point $\varphi =
(\varphi_1,\varphi_2) \in [0,\infty)\! \times\! [0,\infty)$ under
the probability measure $\PP_{\!\varphi}^\infty$ so that the optimal
stopping problem \eqref{3.2} extends as follows
\begin{equation} \h{4pc} \label{3.11}
\hat V(\varphi) = \inf_\tau\;\! \EE_\varphi^\infty \Big[ \int
_0^\tau e^{-\lambda t} \Big( p_1\:\! \varPhi_t^1\! + p_2\:\!
\varPhi_t^2 - \frac{\lambda} {c}\;\! \Big)\;\! dt\:\! \Big]
\end{equation}
for $\varphi \in [0,\infty)\! \times\! [0,\infty)$ with
$\PP_{\!\varphi} (\varPhi_0\! =\! \varphi)=1$ where the infimum is
taken over all stopping times $\tau$ of $\varPhi$ and we recall that
$p_1,p_2 \in [0,1]$ with $p_1 \p p_2 = 1$ are given and fixed. In
this way we have reduced the initial quickest detection problem
\eqref{2.5} or \eqref{2.9} to the optimal stopping problem
\eqref{3.11} for the strong Markov process $\varPhi =
(\varPhi^1,\varPhi^2)$ solving \eqref{3.9}-\eqref{3.10} and being
explicitly given by the Markovian flow \eqref{2.14}+\eqref{2.15} of
the initial point $(\varPhi_0^1,\varPhi_0^2) = (\varphi_1,\varphi_2)
=: \varphi$ in $[0,\infty)\! \times\! [0,\infty)$ under
$\PP_\varphi^\infty$. Note that the optimal stopping problem
\eqref{3.11} is inherently/fully two-dimensional and the
infinitesimal generator of $\varPhi = (\varPhi^1,\varPhi^2)$ is of
elliptic type as discussed in the next section.

\section{Mayer formulation}

The optimal stopping problem \eqref{3.11} is Lagrange formulated. In
this section we derive its Mayer reformulation which is helpful in
the subsequent analysis.

\vspace{4mm}

1.\ From \eqref{3.9}+\eqref{3.10} we read that the infinitesimal
generator of the strong Markov process $\varPhi =
(\varPhi^1,\varPhi^2)$ is given by
\begin{equation} \h{2pc} \label{4.1}
\LL_\varPhi = \lambda\:\! (1 \p \varphi_1)\;\! \partial_{\varphi_1}
+ \lambda\:\! (1 \p \varphi_2)\;\! \partial_{\varphi_2} + \frac{
\mu^2}{2}\:\! \varphi_1^2\, \partial_{\varphi_1 \varphi_1} +
\frac{\mu^2}{2}\:\! \varphi_2^2\, \partial_{\varphi_2 \varphi_2}
\end{equation}
for $(\varphi_1,\varphi_2)$ belonging to $(0,\infty)\! \times\!
(0,\infty)$. From \eqref{2.15} we see that the topological boundary
$\{ 0 \}\! \times\! [0,\infty) \cup (0,\infty)\! \times\! \{ 0 \}$
of the state space $[0,\infty)\! \times\! [0,\infty)$ consists of
\emph{natural} boundary points for $\varPhi$ (meaning that $\varPhi$
can be started at any boundary point never to return to the
boundary) and clearly the differential operator $\LL_\varPhi$ is of
\emph{elliptic} type (cf.\ (2.12) in \cite{Pe-4}).

For the Mayer reformulation of the problem \eqref{3.11} we need to
look for a function $M : [0,\infty)\! \times\! [0,\infty)
\rightarrow \R$ solving the partial differential equation
\begin{equation} \h{9pc} \label{4.2}
\LL_\varPhi M \m \lambda M = L
\end{equation}
on $(0,\infty)\! \times\! (0,\infty)$ where in view of \eqref{3.11}
we set
\begin{equation} \h{7pc} \label{4.3}
L(\varphi_1,\varphi_2) = p_1\:\! \varphi_1 + p_2\:\! \varphi_2 -
\lambda/c
\end{equation}
for $(\varphi_1,\varphi_2) \in [0,\infty)\! \times\! [0,\infty)$.
Ignoring the constant $-\lambda/c$ on the right-hand side of
\eqref{4.2} for now, we see that a possible attempt to solve the
resulting partial differential equation is to separate the variables
$\varphi_1$ and $\varphi_2$ by considering the two ordinary
differential equations
\begin{equation} \h{5pc} \label{4.4}
\lambda\:\! (1 \p \varphi_i)\:\! M_i' + \frac{\mu^2}{2}\:\!
\varphi_i^2\: M_i'' - \lambda\:\! M_i = p_i\:\! \varphi_i
\end{equation}
where $M_i = M_i(\varphi_i)$ is a function/solution to be found for
$\varphi_i \in (0,\infty)$ with $i=1,2$. Simplifying the notation we
see that the equation \eqref{4.4} reads
\begin{equation} \h{7pc} \label{4.5}
x^2 y'' + \kappa\:\! (1 \p x)\:\! y' - \kappa\:\! y = \nu\:\! x
\end{equation}
for $x \in (0,\infty)$ where $y = y(x)$ and we set $\kappa := 2
\lambda/\mu^2$ and $\nu := 2 p_i/\mu^2$ with $i=1,2$. The
homogeneous part of the equation \eqref{4.5} is closely related to
the Euler equation (cf.\ Eq.\ (118) in \cite[Section 2.1.2]{PZ}),
and there exists a general transformation which reduces this part to
another second-order ordinary differential equation, whose leading
term is no longer quadratic but linear, and whose solutions can be
expressed in terms of known special functions (see the reduction of
Eq.\ (129) to Eq.\ (103) and Table 2.2 in \cite[Section 2.1.2]{PZ}).

Motivated by a probabilistic meaning of the posterior probability
distribution ratio process in this context, and aiming to exploit
the specific form of the coefficients $\kappa\:\! (1 \p x)$ and
$-\kappa$ in \eqref{4.5} more directly, we will take a different
tack and seek a solution to \eqref{4.5} by setting
\begin{equation} \h{7pc} \label{4.6}
y(x) := (1 \p x)\, z \Big( \frac{x}{1 \p x} \Big)
\end{equation}
for $x \in(0,\infty)$. Setting $u = x/(1 \p x)$ we then find by
\eqref{4.5} that $z = z(u)$ solves
\begin{equation} \h{7pc} \label{4.7}
u^2 (1 \m u)\:\! z'' + \kappa\:\! z' = \nu\, \frac{u}{1 \m u}
\end{equation}
for $u \in (0,1)$ where the term  $z$ is no longer present. This
equation can therefore be solved in closed form by reduction to a
first-order ordinary differential equation. Inserting this solution
back into \eqref{4.6} we find that the sought solution to
\eqref{4.5} is given by
\begin{equation} \h{2pc} \label{4.8}
y(x) = \nu\;\! (1 \p x)\! \int_0^{x/(1+x)}\!\! \Big( \frac{1
\m v}{v} \Big)^{\!\kappa} e^{\kappa/v}\! \int_0^v \frac{u^{\kappa-1}}
{(1 \m u)^{\kappa+2}}\, e^{-\kappa/u}\, du\;\! dv
\end{equation}
for $x \in (0,\infty)$. This solution can now be used to specify the
sought solutions to the equation \eqref{4.4}. These solutions in
turn can be used to specify the solution to the equation \eqref{4.2}
above. The only matter remaining is to account for the missing
constant $-\lambda/c$ on the right-hand side of \eqref{4.2} and this
will be done shortly below.

\vspace{4mm}

2.\ We now consider the Mayer reformulation of the optimal stopping
problem \eqref{3.11}. Motivated by \eqref{4.8} and recalling that
$\nu = 2 p_i/\mu^2$ with $i=1,2$, let us define a function $M :$
$[0,\infty) \rightarrow \R$ by setting
\begin{equation} \h{2pc} \label{4.9}
M(\varphi) = \frac{2}{\mu^2}\, (1 \p \varphi)\! \int_0^{\varphi/(1
+ \varphi)}\!\! \Big( \frac{1 \m v}{v} \Big)^{\!\kappa} e^{\kappa/v}
\! \int_0^v \frac{u^{\kappa-1}}{(1 \m u)^{\kappa+2}}\, e^{-\kappa/u}\,
du\;\! dv
\end{equation}
for $\varphi \in [0,\infty)$ where we recall that $\kappa = 2
\lambda/\mu^2$. In addition, let us define a function $M :
[0,\infty)\! \times\! [0,\infty) \rightarrow \R$ by setting
\begin{equation} \h{6pc} \label{4.10}
M(\varphi_1,\varphi_2) = p_1\:\! M(\varphi_1) + p_2\:\! M(\varphi_2)
+ 1/c
\end{equation}
for $(\varphi_1,\varphi_2) \in [0,\infty)\! \times\! [0,\infty)$.
The arguments above then show that the function $M$ from
\eqref{4.10} solves the equation \eqref{4.2} above (notice that the
final term $1/c$ yields the missing constant $-\lambda/c$ on the
right-hand side of \eqref{4.2} as needed). Note that we use the same
letter $M$ to denote both functions in order to emphasise the
`fractal' nature of \eqref{4.10} expressed in terms of \eqref{4.9}.
The fact that the two functions have different domains can/will be
used to remove any ambiguity when needed. Having defined the
function $M$ in \eqref{4.10} using \eqref{4.9} we can now describe
the Mayer reformulation of the optimal stopping problem \eqref{3.11}
as follows.

\vspace{6mm}

\textbf{Proposition 2.} \emph{The value function $\hat V$ from
\eqref{3.11} can be expressed as
\begin{equation} \h{6pc} \label{4.11}
\hat V(\varphi) = \inf_\tau\:\! \EE_\varphi^\infty \big[ e^{-\lambda
\tau} M(\varPhi_\tau^1,\varPhi_\tau^2) \big] - M(\varphi)
\end{equation}
for $\varphi \in [0,\infty)\! \times\! [0,\infty)$ where the infimum
is taken over all stopping times $\tau$ of $\varPhi =
(\varPhi^1,\varPhi^2)$ and the function $M$ is given by \eqref{4.10}
using \eqref{4.9} above.}

\vspace{6mm}

\textbf{Proof.} By It\^o's formula using \eqref{3.9}+\eqref{3.10} we
get
\begin{equation} \h{4pc} \label{4.12}
e^{-\lambda t} M(\varPhi_t) = M(\varphi) + \int_0^t e^{-\lambda s}
\big( \LL_\varPhi \m \lambda M)(\varPhi_s)\, ds + N_t
\end{equation}
for $\varphi \in [0,\infty)\! \times\! [0,\infty)$ where $N_t =
\sum_{i=1}^2 \int_0^t e^{-\lambda s} M_{\varphi_i}(\varPhi_s)\:\!
\mu \;\! \varPhi_s^i\, dB_s^i$ is a continuous local martingale for
$t \ge 0$. Making use of a localisation sequence of stopping times
for this local martingale if needed, applying the optional sampling
theorem and recalling that $M$ solves \eqref{4.2}, we find by taking
$\EE_\varphi^\infty$ on both sides in \eqref{4.12} that
\begin{equation} \h{4pc} \label{4.13}
\EE_\varphi^\infty \big[ e^{-\lambda \tau} M(\varPhi_\tau^1,
\varPhi_\tau^2) \big] = M(\varphi) + \EE_\varphi^\infty \Big[
\int_0^\tau\!\! e^{-\lambda t} L(\varPhi_t)\, dt\:\! \Big]
\end{equation}
for all $\varphi \in [0,\infty)\! \times\! [0,\infty)$ and all
(bounded) stopping times $\tau$ of $\varPhi$. From \eqref{3.11} and
\eqref{4.13} using \eqref{4.3} we see that \eqref{4.11} holds as
claimed and the proof is complete. \hfill $\square$

\vspace{6mm}

3.\ From Proposition 2 we see that the optimal stopping problem
\eqref{3.11} is equivalent to the optimal stopping problem defined
by
\begin{equation} \h{6pc} \label{4.14}
\check V(\varphi) = \inf_\tau\:\! \EE_\varphi^\infty \big[ e^{-\lambda
\tau} M(\varPhi_\tau^1,\varPhi_\tau^2) \big]
\end{equation}
for $\varphi \in [0,\infty)\! \times\! [0,\infty)$ where the infimum
is taken over all stopping times $\tau$ of $\varPhi =
(\varPhi^1,\varPhi^2)$ and the function $M$ is given by \eqref{4.10}
using \eqref{4.9} above. The optimal stopping problem \eqref{4.14}
is Mayer formulated. From \eqref{4.11} and \eqref{4.14} we see that
\begin{equation} \h{7.5pc} \label{4.15}
\hat V(\varphi) = \check V(\varphi) - M(\varphi)
\end{equation}
for $\varphi \in [0,\infty)\! \times\! [0,\infty)$. The Mayer
reformulation \eqref{4.14} has certain advantages that will be
exploited in the subsequent analysis of the optimal stopping problem
\eqref{3.11} below.

\section{One dimension}

The observed process $X$ in the initial quickest detection problem
\eqref{2.5} is two-dimensional. In this section we consider the
analogue of \eqref{2.5} and the resulting optimal stopping problem
\eqref{3.11} when $X$ is one-dimensional. The reduction of dimension
from two to one corresponds to taking either $p_1$ or $p_2$ equal to
$1$. Then $\varPhi$ standing for either $\varPhi^1$ or $\varPhi^2$
respectively is a one-dimensional Markov/diffusion process so that
standard optimal stopping arguments can be used to solve the
problem. The derived results for the one-dimensional optimal
stopping problem \eqref{3.11} when $X$ is one-dimensional will be
used in the subsequent analysis of the two-dimensional optimal
stopping problem \eqref{3.11} when $X$ is two-dimensional.

\vspace{4mm}

1.\ Using the same arguments as in Sections 2 and 3 above, it is
easily seen that the quickest detection problem \eqref{2.5} when $X$
is one-dimensional reduces to the optimal stopping problem
\eqref{3.11} with $p_1=1$ and $p_2=0$ (without loss of generality).
Omitting the superscript $1$ from $\varPhi^1$ for simplicity, we
thus see that the optimal stopping problem \eqref{3.11} reads
\begin{equation} \h{5pc} \label{5.1}
\hat V(\varphi) = \inf_\tau\;\! \EE_\varphi^\infty \Big[ \int
_0^\tau e^{-\lambda t} \Big( \varPhi_t - \frac{\lambda}{c}
\;\! \Big)\;\! dt\:\! \Big]
\end{equation}
for $\varphi \in [0,\infty)$ with $\PP_\varphi^\infty(\varPhi_0\!
=\! \varphi)=1$ where the infimum is taken over all stopping times
$\tau$ of $\varPhi$. From \eqref{3.9} we see that the infinitesimal
generator of $\varPhi$ is given by
\begin{equation} \h{5pc} \label{5.2}
\LL_\varPhi = \lambda\:\! (1 \p \varphi)\;\! \frac{d}{d \varphi}
+ \frac{\mu^2}{2}\:\! \varphi^2\, \frac{d^2}{d \varphi^2}
\end{equation}
for $\varphi$ belonging to $[0,\infty)$.

\vspace{4mm}

2.\ Noting that the optimal stopping problem \eqref{5.1} is Lagrange
formulated, standard arguments imply (see e.g.\ \cite{PS}) that
$\hat V$ should solve the free-boundary problem
\begin{align} \h{5pc} \label{5.3}
&\LL_\varphi \hat V \m \lambda \hat V = -(\varphi \m \lambda/c)\;\;
\text{for}\;\; \varphi \in [0,\varphi_*) \\[1pt] \label{5.4} &\hat V(
\varphi_*) = 0\;\; \text{(instantaneous stopping)} \\[1pt] \label{5.5}
&\hat V'(\varphi_*) = 0\;\; \text{smooth fit}
\end{align}
where $\varphi_* \in (\lambda/c,\infty)$ is the optimal stopping
boundary/point to be found, and we set $\hat V(\varphi)=0$ for
$\varphi \in (\varphi_*,\infty)$ in addition to \eqref{5.4} above
(note from \eqref{5.1} that considering the exit times of $\varPhi$
from sufficiently small intervals shows that it is never optimal to
stop at least in $[0,\lambda/c)$ as partly indicated above).

\vspace{4mm}

3.\ The general solution to the ordinary differential equation
\eqref{5.3} is given by
\begin{equation} \h{2pc} \label{5.6}
\hat V(\varphi) = A\;\! (1 \p \varphi) \int_{1/2}^{\varphi/(1+\varphi)}
\!\! \Big( \frac{1 \m v}{v} \Big)^{\!\kappa} e^{\kappa/v}\:\! dv +
B\;\! (1 \p \varphi) + \hat V_p(\varphi)
\end{equation}
where $A$ and $B$ are (unspecified) real constants and $V_p$ is a
particular solution to \eqref{5.3} obtained by modifying the
solution $M$ from \eqref{4.9} as follows
\begin{equation} \h{2pc} \label{5.7}
\hat V_p(\varphi) = -\frac{2}{\mu^2}\, (1 \p \varphi)\! \int_{
\varphi_* /(1 + \varphi_*)}^{\varphi/(1 + \varphi)}\!\! \Big(
\frac{1 \m v}{v} \Big)^{\!\kappa} e^{\kappa/v} \! \int_0^v \frac{
u^{\kappa-1}} {(1 \m u)^{\kappa+2}}\, e^{-\kappa/u}\, du\;\! dv - 1/c
\end{equation}
for $\varphi \in [0,\varphi_*]$ where $\varphi_* \in
(\lambda/c,\infty)$ is to be found and we recall that $\kappa = 2
\lambda/\mu^2$. This can be obtained by noticing that $M$ from
\eqref{4.10} with $p_1=1$ and $p_2=0$ solves \eqref{4.2} with
$\LL_\varPhi$ from \eqref{5.2} if and only if $\hat V := -M$ solves
\eqref{5.3}. Hence we see that the transformation \eqref{4.6}
reduces the equation \eqref{5.3} to a solvable form. Proceeding thus
as in \eqref{4.7} and \eqref{4.8} above, and applying the analogous
arguments to the homogeneous part of the equation \eqref{5.3}, we
obtain the general solution \eqref{5.6} with \eqref{5.7} as claimed.
The motivation for modifying the particular solution $M$ from
\eqref{4.9} as $\hat V_p$ in \eqref{5.7} comes from the
instantaneous stopping and smooth fit conditions \eqref{5.4} and
\eqref{5.5} as will be clear from the calculations below.

\vspace{2mm}

4.\ A direct differentiation in \eqref{5.6} shows that $\hat V'(0+)
= +\infty$ if $A>0$ and $\hat V'(0+) = -\infty$ if $A<0$. We thus
choose $A=0$ as a candidate value in the sequel. Using \eqref{5.4}
we then find that $B=1/c(1 \p \varphi_*)$. This yields the following
candidate solution to the free-boundary problem
\eqref{5.3}-\eqref{5.5} above
\begin{equation} \h{0pc} \label{5.8}
\hat V(\varphi) = - \frac{\varphi_* \m \varphi}{c(1 \p \varphi_*)}
- \frac{2}{\mu^2}\, (1 \p \varphi)\! \int_{\varphi_*/(1 + \varphi_*)}
^{\varphi/(1 + \varphi)}\!\! \Big(\frac{1 \m v}{v} \Big)^{\!\kappa}
e^{\kappa/v} \! \int_0^v \frac{ u^{\kappa-1}} {(1 \m u)^{\kappa+2}}
\, e^{-\kappa/u}\, du\;\! dv
\end{equation}
for $\varphi \in [0,\varphi_*]$ and $\hat V(\varphi)=0$ for
$(\varphi_*,\infty)$. A direct differentiation in \eqref{5.8} then
shows that \eqref{5.5} holds if and only if $\varphi_*$ solves
\begin{equation} \h{5pc} \label{5.9}
\frac{e^{\kappa(1+\varphi_*)/\varphi_*}}{\varphi_*^\kappa} \int_0^{
\varphi_*/(1 + \varphi_*)}\!\! \frac{u^{\kappa-1}} {(1 \m u)^{\kappa+2}}
\, e^{-\kappa/u}\, du = \frac{\mu^2}{2c}
\end{equation}
where we recall that $\kappa = 2 \lambda/\mu^2$.

\vspace{4mm}

5.\ To make the arguments developed above rigorous we can reverse
their order and start our analysis from the end. Firstly, we claim
that there exists a unique point $\varphi_* \in (\lambda/c,\infty)$
satisfying the equation \eqref{5.9}. For this, define the functions
$F$ and $G$ by setting
\begin{equation} \h{0pc} \label{5.10}
F(\varphi) = \int_0^{\varphi/(1+\varphi)}\!\! \frac{1}{u(1 \m u)^2}
\, e^{\kappa(\log[u/(1-u)]-1/u)}\;\! du\;\;\; \&\;\;\; G(\varphi) =
\frac{\mu^2}{2c}\, e^{\kappa(\log \varphi - (1 + \varphi)/\varphi)}
\end{equation}
for $\varphi \in [0,\infty)$. Note that the claim about \eqref{5.9}
is equivalent to establishing that
\begin{equation} \h{9pc} \label{5.11}
F(\varphi_*) = G(\varphi_*)
\end{equation}
for a unique point $\varphi_* \in (\lambda/c,\infty)$. To verify
\eqref{5.11} note that $F(0)=G(0)$ and we have
\begin{equation} \h{3pc} \label{5.12}
\varphi < \frac{\lambda}{c} \Longleftrightarrow F'(\varphi)
< G' (\varphi)\;\;\; \&\;\;\; \varphi > \frac{\lambda}{c}
\Longleftrightarrow F'(\varphi) > G' (\varphi)
\end{equation}
for $\varphi \in (0,\infty)$ as is easily verified by a direct
differentiation in \eqref{5.10}. This shows that $F(\varphi) <
G(\varphi)$ for all $\varphi \in (0,\lambda/c]$. Moreover, applying
L'Hospital's rule we find that
\begin{equation} \h{5pc} \label{5.13}
\lim_{\varphi \rightarrow \infty} \frac{F(\varphi)}{G(\varphi)} =
\lim_{\varphi \rightarrow \infty} \frac{F'(\varphi)}{G'(\varphi)}
= \lim_{\varphi \rightarrow \infty} \Big( \frac{c}{\lambda}\;\!
\varphi \Big) = \infty\, .
\end{equation}
Combining \eqref{5.12} and \eqref{5.13} we see that the graphs of
$F$ and $G$ must intersect on $(\lambda/c,\infty)$ at a unique point
$\varphi_*$ establishing \eqref{5.11} as claimed. Secondly, define
$\hat V_*(\varphi)$ by the right-hand side of \eqref{5.8} for
$\varphi \in [0,\varphi_*]$ and set $\hat V_*(\varphi)=0$ for
$(\varphi_*,\infty)$. Then the arguments above show (or it is a
matter of routine to verify) that $\hat V_*$ solves the
free-boundary problem \eqref{5.3}-\eqref{5.5} above. Thirdly,
applying the It\^o-Tanaka formula (cf.\ \cite[p.\ 223]{RY}) to $\hat
V_*$ composed with $\varPhi$, which reduces to It\^o's formula due
to smooth fit \eqref{5.5}, and making use of the optional sampling
theorem, it is easily verified that $\hat V_*$ coincides with the
value function $\hat V$ from \eqref{5.1} and the optimal stopping
time (at which the infimum in \eqref{5.1} is attained) is given by
\begin{equation} \h{7pc} \label{5.14}
\tau_* = \inf\, \{\, t \ge 0\, \vert\, \varPhi_t \in [\varphi_*,
\infty)\, \}
\end{equation}
where $\varphi_* \in (\lambda/c,\infty)$ is a unique solution to
\eqref{5.9} on $(0,\infty)$. These facts will be used in the
subsequent analysis of the optimal stopping problem \eqref{3.11}
when the observed process $X$ is two-dimensional as assumed in
Section 2 above.

\section{Properties of the optimal stopping boundary}

In this section we establish the existence of an optimal stopping
time in the problem \eqref{3.11} and derive basic properties of the
optimal stopping boundary.

\vspace{2mm}

1.\ Looking at \eqref{3.11} we may conclude that the (candidate)
continuation and stopping sets in this problem need to be defined as
follows
\begin{align} \h{5.5pc} \label{6.1}
C = \{\, \varphi \in [0,\infty)\! \times\! [0,\infty)\; \vert\;
\hat V(\varphi)<0\, \} \\[2pt] \label{6.2} D = \{\, \varphi \in
[0,\infty)\! \times\! [0,\infty)\; \vert\; \hat V(\varphi)<0\, \}
\end{align}
respectively. Recalling that \eqref{2.15} defines a Markovian
functional of the initial point $\varPhi_0^i := \varphi_i$ in
$[0,\infty)$ of the process $\varPhi^i$ for $i=1,2$, we see that the
expectation in \eqref{4.11} defines a continuous function of the
initial point $\varphi = (\varphi_1,\varphi_2)$ of the process
$\varPhi = (\varPhi^1,\varPhi^2)$ for every (bounded) stopping time
$\tau$ of $\varPhi$ given and fixed. Taking the infimum over all
(bounded) stopping times $\tau$ of $\varPhi$ we can thus conclude
from \eqref{4.11} that the value function $\hat V$ is upper
semicontinuous on $[0,\infty)\! \times\! [0,\infty)$. From
\eqref{4.10} with \eqref{4.9} we see that the loss function $M$ in
\eqref{4.11} is continuous and hence lower semicontinuous too. It
follows therefore by \cite[Corollary 2.9]{PS} that the first entry
time of the process $\varPhi$ into the closed set $D$ defined by
\begin{equation} \h{7pc} \label{6.3}
\tau_D = \inf\, \{\, t \ge 0\; \vert\; \varPhi_t \in D\, \}
\end{equation}
is optimal in \eqref{4.11} and hence in \eqref{3.11} as well
whenever $\PP_{\!\varphi}( \tau_D < \infty ) = 1$ for $\varphi \in
[0,\infty)\! \times [0,\infty)$. In the sequel we will establish
this and other properties of $\tau_D$ by analysing the boundary of
$D$.

\vspace{2mm}

2.\ To derive an upper bound on the boundary of $D$, recall that the
optimal stopping boundary/point $\varphi_*$ in the one-dimensional
problem \eqref{5.1} can be characterised as a unique solution to
\eqref{5.9} on $(\lambda/c,\infty)$. Note from \eqref{5.9} that
$\varphi_* = \varphi_*(\lambda/\mu^2,\lambda/c)$ and set
\begin{equation} \h{5.5pc} \label{6.4}
\varphi_1^* := \varphi_* \big(\tfrac{\lambda}{\mu^2},\tfrac{\lambda}
{p_1 c} \big)\;\;\; \&\;\;\; \varphi_2^* := \varphi_* \big(\tfrac
{\lambda}{\mu^2},\tfrac{\lambda} {p_2 c} \big)
\end{equation}
for $\lambda>0$, $\mu \in \R$,$c>0$ and $p_1, p_2 \in (0,1)$ with
$p_1 \p p_2 = 1$. Recall that $\varphi_1^* \in (\lambda/(p_1
c),\infty)$ and $\varphi_2^* \in (\lambda/(p_2 c),\infty)$. We can
now expose basic properties of the the value function and the
continuation/stopping set in the problem \eqref{3.11} as follows. \\\\
\textbf{Proposition 3.} 
\begin{align}\label{6.5} \textit{The value function}\;\; \hat V\;\; \textit{is concave
and continuous on}\;\; [0,\infty)\! \times\! [0,\infty)\, .
\end{align}

\begin{figure}[H] 
\begin{center}
\includegraphics[scale=0.8]{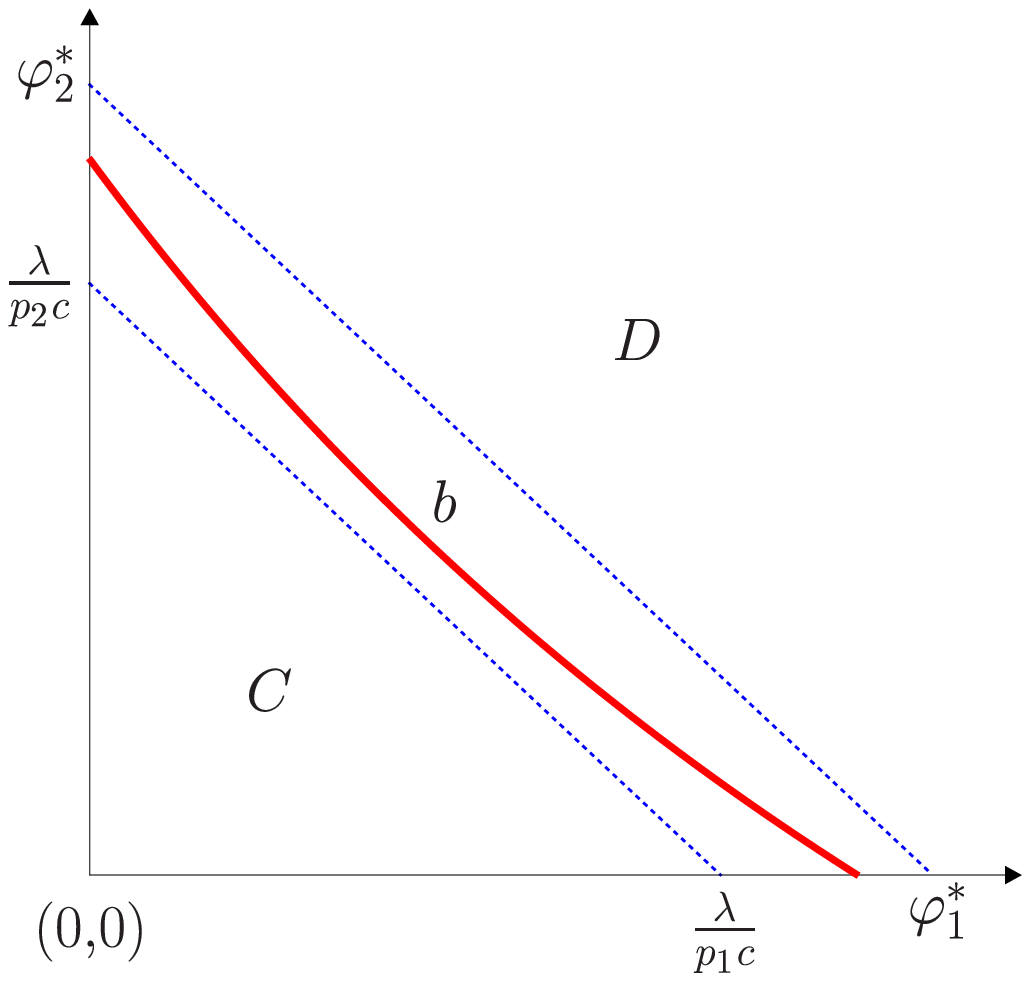}
\end{center}

{\par \leftskip=2cm \rightskip=2cm \small \noindent \vspace{-2mm}

\textbf{Figure 1.} The optimal stopping boundary $b$ in the problem
\eqref{3.11} when $\mu=\lambda=c=1$ and $p_1=p_2=1/2$. The dotted
lines are graphs of the linear functions from \eqref{6.11} below.

\par} \vspace{0mm}

\end{figure}
\emph{\begin{align}
\label{6.6} &\h{-5pt}\text{If}\;\; \varphi_1 \le \psi_1\;\; \&\;\;
\varphi_2 \le \psi_2\;\; \text{then}\;\; \hat V(\varphi_1,\psi_1)
\le \hat V(\varphi_2,\psi_2)\, . \\[3pt] \label{6.7} &\h{-5pt}
\text{If}\;\; (\varphi_1,\varphi_2) \in D\;\; \text{and}\;\; \psi_1
\ge \varphi_1\;\; \&\;\; \psi_2 \ge \varphi_2\;\; \text{then}\;\;
(\psi_1,\psi_2) \in D\, . \\[3pt] \label{6.8} &\h{-5pt}\text{The
stopping set}\;\; D\;\; \text{is convex}\;\; \text{and the
trigon}\;\; \{\, (\varphi_1,\varphi_2) \in [0,\infty)\! \times\!
[0,\infty)\; \vert\; \\[-3pt] \notag &\h{-3pt}\varphi_1/\varphi_1^*
\p \varphi_2/\varphi_2^* \m 1 \ge 0\, \} \;\; \text{is contained
in}\;\; D\, . \\[3pt] \label{6.9} &\h{-5pt} \text{The triangle}
\;\; \{\, (\varphi_1,\varphi_2) \in [0,\infty) \! \times\! [0,
\infty)\; \vert\; p_1 \varphi_1 \p p_2 \varphi_2 \m \lambda/c <
0\, \}\;\; \text{is contained} \\[-3pt] \notag &\h{-5pt}\text{in
the continuation set}\;\; C\, .
\end{align}}


\vspace{-4mm}

\textbf{Proof.} \eqref{6.5}: Combining the fact that the Markovian
flow \eqref{2.15} is linear as a function of its initial point with
the fact that the integral in \eqref{3.11} is a linear function of
its argument, and using that the infimum of a convex combination is
larger than the convex combination of the infima, we find that $\hat
V$ is concave on $[0,\infty)\! \times\! [0,\infty)$ as claimed.
Hence we can also conclude that $\hat V$ is continuous on the open
set $(0,\infty)\! \times\! (0,\infty)$. To see that $\hat V$ is
continuous at the boundary points of $[0,\infty)\! \times\!
[0,\infty)$ we may recall the well-known (and easily verified) fact
that the concave function $\hat V$ is lower semicontinuous on the
closed and convex set $[0,\infty)\! \times\! [0,\infty)$. Since we
also know that $\hat V$ is upper semicontinuous on $[0,\infty)\!
\times\! [0,\infty)$ as established following \eqref{6.2} above, we
see that $\hat V$ is continuous on the entire $[0,\infty)\! \times\!
[0,\infty)$ as claimed.

\vspace{2mm}

\eqref{6.6}: This is a direct consequence of the fact that the
Markovian flow \eqref{2.15} is increasing as a function of its
initial point being used in \eqref{3.11} above.

\vspace{4mm}

\eqref{6.7}: By \eqref{6.6} we have $\hat V(\varphi_1,\varphi_2) \le
\hat V(\psi_1,\psi_2) \le 0$ so that $(\varphi_1,\varphi_2) \in D$
i.e.\ $\hat V(\varphi_1,\varphi_2) = 0$ implies that $\hat
V(\psi_1,\psi_2) = 0$ i.e.\ $(\psi_1,\psi_2) \in D$ as claimed.

\vspace{4mm}

\eqref{6.8}: To see that $D$ is convex, take any $\varphi$ and
$\psi$ from $D$ and note by \eqref{6.5} that $0 \ge$ $\hat V(\alpha
\varphi \p (1 \m \alpha) \psi) \ge \alpha \hat V(\varphi) \p (1 \m
\alpha) \hat V(\psi) = 0$ so that $\hat V(\alpha \varphi \p (1 \m
\alpha) \psi) = 0$ i.e.\ $\alpha \varphi \p (1 \m \alpha) \psi \in
D$ for every $\alpha \in [0,1]$ as claimed. To see that the trigon
is contained in $D$, note that pulling $p_1$ in front of the infimum
in \eqref{3.11} shows that the point $(\varphi_1^*,0)$ belongs to
$D$ because $\varphi_1^*$ as defined in \eqref{6.4} above is an
optimal stopping point in the one-dimensional problem obtained by
removing the (independent) positive term $(p_2/p_1)\;\! \varPhi_t^2$
from the integral with respect to time in \eqref{3.11} with $p_1$ in
front of the infimum. Similarly, we see that the point
$(0,\varphi_2^*)$ belongs to $D$. But then the entire trigon is
contained in $D$ due to its convexity.

\vspace{4mm}

\eqref{6.9}: Taking any point $\varphi$ from the triangle and
replacing $\tau$ in \eqref{3.11} by the first exit time of $\varPhi$
from a sufficiently small ball around $\varphi$ that is strictly
contained in the triangle, we see that the integrand in \eqref{3.11}
remains strictly negative so that $\hat V$ takes a strictly negative
value at $\varphi$ itself, showing that $\varphi$ belongs to the
continuation set $C$ as claimed. \hfill $\square$

\vspace{4mm}

3.\ From the results of Proposition 3 we see that the stopping set
in the problem \eqref{3.11} can be described  as  follows
\begin{equation} \h{4.5pc} \label{6.10}
D = \{\, (\varphi_1,\varphi_2) \in [0,\infty)\! \times\! [0,\infty)
\; \vert\; \varphi_2 \ge b(\varphi_1)\, \}
\end{equation}
where $b : [0,\infty) \rightarrow [0,\infty)$ is a convex,
continuous, decreasing function satisfying
\begin{equation} \h{5.5pc} \label{6.11}
-\frac{p_1}{p_2}\, \varphi_1 + \frac{\lambda}{p_2 c} \le b(\varphi_1)
\le - \frac{\varphi_2^*}{\varphi_1^*}\, \varphi_1 + \varphi_2^*
\end{equation}
for $\varphi_1 \in [0,\lambda/(p_1 c)]$ in the first inequality and
$\varphi_1 \in [0,\varphi_1^*]$ in the second inequality
respectively (see Figure 1). Note that since $\varPhi^1$ and
$\varPhi^2$ are independent and either recurrent or transient in
$[0,\infty)$ (converging to $\infty$ in the latter case) as recalled
following \eqref{3.9}+\eqref{3.10} above, we see from \eqref{6.11}
that $\PP_{\!\varphi}( \tau_D < \infty ) = 1$ for all $\varphi \in
[0,\infty)\! \times [0,\infty)$ as claimed following \eqref{6.3}
above. We address the question of characterising/determining $b$ in
the remaining two sections.

\section{Free-boundary problem}

In this section we derive a free-boundary problem that stands in
one-to-one correspondence with the optimal stopping problem
\eqref{3.11}. Using the results derived in the previous sections we
show that the value function $\hat V$ from \eqref{3.11} and the
optimal stopping boundary $b$ from \eqref{6.10} solve the
free-boundary problem. This establishes the existence of a solution
to the free-boundary problem. Its uniqueness in a natural class of
functions will follow from a more general uniqueness result that
will be established in Section 8 below. This will also yield an
explicit  integral representation of the value function $\hat V$
expressed in terms of the optimal stopping boundary $b$.

\vspace{2mm}

1.\ Consider the optimal stopping problem \eqref{3.11} where the
Markov process $\varPhi = (\varPhi^1,\varPhi^2)$ solves the system
of stochastic differential equations \eqref{3.9}-\eqref{3.10} driven
by a standard Brownian motion $B = (B^1,B^2)$ under the probability
measure $\PP^\infty$. Recall that the infinitesimal generator of
$\varPhi$ is the second-order elliptic differential operator
$\LL_\varPhi$ given in \eqref{4.1} above. Looking at \eqref{3.11}
and relying on other properties of $\hat V$ and $b$ derived above,
we are naturally led to formulate the following free-boundary
problem for finding $\hat V$ and $b$ :
\begin{align} \h{4pc} \label{7.1}
&\LL_\varPhi \hat V \m \lambda \hat V = -L\;\;\; \text{in}\;\; C
\\ \label{7.2} &\hat V(\varphi) = 0\;\;\; \text{for}\;\; \varphi
\in D\;\;\; \text{(instantaneous stopping)} \\ \label{7.3} &\hat
V_{\varphi_i}(\varphi) = 0\;\;\; \text{for}\;\; \varphi \in
\partial C\;\; \text{and} \;\; i=1,2\;\;\; \text{(smooth fit)}
\end{align}
where $L$ is defined in \eqref{4.3} above, $C$ is the (continuation)
set from \eqref{6.1} above, $D$ is the (stopping) set from
\eqref{6.2}+\eqref{6.10} above, and $\partial C = \{\,
(\varphi_1,\varphi_2) \in [0,\infty)\! \times\! [0,\infty)\; \vert\;
\varphi_2 = b(\varphi_1)\, \}$ is the (optimal stopping) boundary
between the sets $C$ and $D$.

\vspace{4mm}

2.\ To formulate the existence and uniqueness result for the
free-boundary problem \eqref{7.1}-\eqref{7.3}, we let $\cal C$
denote the class of functions $(U,a)$ such that
\begin{align} \h{0pc} \label{7.4}
&U\;\; \text{belongs to}\;\; C^1(\bar C_a) \cap C^2(C_a)\;\;
\text{and}\;\; \text{is continuous \& bounded on}\;\; [0,\infty)
\! \times\! [0,\infty) \\[3pt] \label{7.5} &b\;\; \text{is continuous
\& decreasing on}\;\; [0,\infty)\;\; \text{and satisfies}\;\; p_1
\varphi_1 \p p_2 b(\varphi_1) \m \lambda/c \ge 0 \\[-3pt] \notag
&\text{for}\;\; \varphi_1 \in [0,\infty)
\end{align}
where we set $C_a = \{\, (\varphi_1,\varphi_2) \in [0,\infty)\!
\times\! [0,\infty)\; \vert\; \varphi_2 < a(\varphi_1)\, \}$ and
$\bar C_a = \{\, (\varphi_1,\varphi_2) \in [0,\infty)\! \times\!
[0,\infty)\; \vert\; \varphi_2 \le a(\varphi_1)\;\; \& \;\;
\varphi_1 \le \inf \{\, \psi_1 \in [0,\infty)\; \vert\; a(\psi_1) =
0\, \}\, \}$. Note that in the latter set we only account for the
smallest zero of the function $a$ should such zeros exist.

\vspace{6mm}

\textbf{Theorem 4.} \emph{The free-boundary problem
\eqref{7.1}-\eqref{7.3} has a unique solution $(\hat V,b)$ in the
class $\cal C$ where $\hat V$ is given in \eqref{3.11} and $b$ is
given in \eqref{6.10} above.}

\vspace{6mm}

\textbf{Proof.} We first show that the pair $(\hat V,b)$ belongs to
the class $\cal C$ and solves the free-boundary problem
\eqref{7.1}-\eqref{7.3}. For this, note that the optimal stopping
problem \eqref{3.11} is Lagrange formulated so that standard
arguments (see e.g.\ the final paragraph of Section 2 in
\cite{DePe}) imply that $\hat V$ belongs to $C^2(C)$ and satisfies
\eqref{7.1}. From \eqref{6.5} we know that $\hat V$ is continuous on
$[0,\infty) \! \times\! [0,\infty)$ and from \eqref{3.11} we readily
find that
\begin{equation} \h{9pc} \label{7.6}
- \frac{1}{c} \le \hat V(\varphi) \le 0
\end{equation}
for all $\varphi \in [0,\infty)\! \times\! [0,\infty)$. Moreover,
recall that the process $\varPhi = (\varPhi^1,\varPhi^2)$ is strong
Feller while it is evident that each point $\varphi \in \partial C$
is probabilistically regular for the set $D$ since $b$ is decreasing
and the coordinate processes $\varPhi^1$ \& $\varPhi^2$ are
independent. Finally, from \eqref{2.15} we see that the process
$\varPhi$ can be realised as a continuously differentiable
stochastic flow of its initial point so that the integrability
conditions of Theorem 8 in \cite{DePe} are satisfied. Recalling that
$\hat V$ satisfies \eqref{7.2}, and applying the result of that
theorem, we can conclude that
\begin{equation} \h{4pc} \label{7.7}
\hat V\;\; \text{is continuously differentiable on}\;\; [0,\infty)
\! \times\! [0,\infty)\, .
\end{equation}
In particular, this shows that \eqref{7.3} holds as well as that
$\hat V$ belongs to $C^1(\bar C)$ as required in \eqref{7.4} above.
The fact that $b$ satisfies \eqref{7.5} was established in the final
paragraph of Section 6 above. This shows that $(\hat V,b)$ belongs
$\cal C$ and solves \eqref{7.1}-\eqref{7.3} as claimed. To derive
uniqueness of the solution we will first see in the next section
that any solution $(U,a)$ to \eqref{7.1}-\eqref{7.3} from the class
$\cal C$ admits an explicit integral representation for $U$
expressed in terms of $a$, which in turn solves a nonlinear Fredholm
integral equation, and we will see that this equation cannot have
other solutions satisfying the required properties. From these facts
we can conclude that the free-boundary problem
\eqref{7.1}-\eqref{7.3} cannot have other solutions in the class
$\cal C$ as claimed. This completes the proof. \hfill $\square$

\section{Nonlinear integral equation}

In this section we show that the optimal stopping boundary $b$ from
\eqref{6.10} can be characterised as the unique solution to a
nonlinear Fredholm integral equation. This also yields an explicit
integral representation of the value function $\hat V$ from
\eqref{3.11} expressed in terms of the optimal stopping boundary
$b$. As a consequence of the existence and uniqueness result for the
the nonlinear Fredholm integral equation we also obtain uniqueness
of the solution to the free-boundary problem \eqref{7.1}-\eqref{7.3}
as explained in the proof of Theorem 4 above. Finally, collecting
the results derived throughout the paper we conclude our exposition
by disclosing the solution to the initial problem.

\vspace{6mm}

1.\ Let $p = p(t;\varphi_1,\varphi_2,\psi_1,\psi_2)$ denote the
transition probability density function of the Markov process
$\varPhi = (\varPhi^1,\varPhi^2)$ in the sense that
\begin{equation} \h{4pc} \label{8.1}
\PP_{\!\!\varphi_1,\varphi_2}^\infty \big( \varPhi_t\! \in\! A \big)
= \iint_A p(t;\varphi_1,\varphi_2,\psi_1,\psi_2)\, d\psi_1 d\psi_2
\end{equation}
for any measurable $A \subseteq [0,\infty)\! \times\! [0,\infty)$
with $(\varphi_1,\varphi_2) \in [0,\infty)\! \times\! [0,\infty)$
and $t \ge 0$ given and fixed. Since $\varPhi^1$ and $\varPhi^2$ are
independent, we have $p(t;\varphi_1,\varphi_2,\psi_1,\psi_2) =
p_1(t;\varphi_1,\psi_1)\, p_2(t;\varphi_2,\psi_2)$ for all
$(\varphi_1,\varphi_2)$ \& $(\psi_1,\psi_2)$ in $[0,\infty)\!
\times\! [0,\infty)$ and all $t \ge 0$, where $p_1$ and $p_2$ are
transition probability density functions of $\varPhi^1$ and
$\varPhi^2$ respectively. Explicit expressions for $p_1$ and $p_2$
are known (see e.g.\ \cite{Pe-2} and the references therein). Having
$p$ we can evaluate the expression of interest in the theorem below
as follows
\begin{align} \h{4pc} \label{8.2}
K(t;\varphi_1,\varphi_2) &:= \EE_{\varphi_1,\varphi_2}^\infty \big[
L(\varPhi_t^1,\varPhi_t^2)\;\! I \big( \varphi_t^2\! <\! b(\varphi_t^1)
\big) \big] \\[3pt] \notag &= \int_0^{\varphi_0}\! d\psi_1\! \int_0^{b(
\psi_1)}\!\! L(\psi_1, \psi_2)\, p(t;\varphi_1,\varphi_2,\psi_1,
\psi_2)\, d\psi_2
\end{align}
for $t \ge 0$ and $(\varphi_1,\varphi_2)\in [0,\infty)\! \times\!
[0,\infty)$ where $\varphi_0$ is the smallest zero of $b$ on
$[0,\infty)$ (recall that $\varphi_0 \in [\lambda/(p_1
c),\varphi_1^*]$ as seen in Figure 1 above) and $L$ is defined in
\eqref{4.3} above.

\vspace{6mm}

\textbf{Theorem 5 (Existence and uniqueness).} \emph{The optimal stopping
boundary $b$ in \eqref{3.11} can be characterised as the unique
solution to the nonlinear Fredholm integral equation
\begin{equation} \h{7pc} \label{8.3}
\int_0^\infty\! e^{-\lambda t} K(t;\varphi_1,b(\varphi_1))\, dt = 0
\end{equation}
in the class of continuous {\rm \&} decreasing (convex) functions
$b$ on $[0,\infty)$ satisfying $p_1 \varphi_1 \p p_2 b(\varphi_1) \m
\lambda/c \ge 0$ for $\varphi_1 \in [0,\varphi_0)$ where $\varphi_0$
is the smallest zero of $b$ on $[0,\infty)$. The value function
$\hat V$ in \eqref{3.11} admits the following representation
\begin{equation} \h{7pc} \label{8.4}
\hat V(\varphi_1,\varphi_2) = \int_0^\infty\! e^{-\lambda t} K(t;
\varphi_1,\varphi_2)\, dt
\end{equation}
for $(\varphi_1,\varphi_2)\in [0,\infty)\! \times\! [0,\infty)$. The
optimal stopping time in \eqref{3.11} is given by
\begin{equation} \h{7pc} \label{8.5}
\tau_b = \inf\, \{\, t \ge 0\; \vert\; \varPhi_t^2 \ge b(\varPhi
_t^1)\, \}
\end{equation}
under $\PP_{\!\!\varphi_1,\varphi_2}^\infty$ with
$(\varphi_1,\varphi_2)\in [0,\infty)\! \times\! [0,\infty)$ given
and fixed.}

\vspace{6mm}

\textbf{Proof.} 1.\ \emph{Existence}. We first show that the optimal
stopping boundary $b$ in \eqref{3.11} solves \eqref{8.3}. Recalling
that $b$ satisfies the properties stated following \eqref{6.10}
above, this will establish the existence of a solution to
\eqref{8.3} in the specified class of functions.

For this, to gain control over the (individual) second partial
derivatives $\hat V_{\varphi_1 \varphi_1}$ and $\hat V_{\varphi_2
\varphi_2}$ close to the optimal stopping boundary within $C$ (see
\cite{GT} for general results of this kind), consider the sets $C_n
:= \{\, \varphi \in [0,\infty)\! \times\! [0,\infty)\; \vert\; \hat
V(\varphi) < -1/n\, \}$ and $D_n := \{\, \varphi \in [0,\infty)\!
\times\! [0,\infty)\; \vert\; \hat V(\varphi) \ge -1/n\, \}$ for $n
\ge 1$ (large). Note that $C_n \uparrow C$ and $D_n \downarrow D$ as
$n \uparrow \infty$. Moreover, using the same arguments as for the
sets $C$ and $D$ above, we find that the set $D_n$ is convex, and
the boundary $b_n = b_n(\varphi_1)$ between $C_n$ and $D_n$ is a
convex, continuous, decreasing function of $\varphi_1$ in
$[0,\varphi_0^n]$ where $\varphi_0^n$ is the smallest zero of $b_n$
on $[0,\varphi_1^*]$ for $n \ge 1$. This also shows that $b_n
\uparrow b$ uniformly on $[0,\varphi_0^n]$ with $\varphi_0^n
\uparrow \varphi_0$ as $n \rightarrow \infty$ where $\varphi_0$ is
the smallest zero of $b$ on $[0,\varphi_1^*]$.

Approximate the value function $\hat V$ in \eqref{3.11} by functions
$\hat V^n$ defined as $\hat V$ on $C_n$ and $-1/n$ on $D_n$ for $n
\ge 1$. Note that $\hat V^n \uparrow \hat V$ uniformly on
$[0,\infty)\! \times\! [0,\infty)$ as $n \rightarrow \infty$.
Moreover, letting $n \ge 1$ be given and fixed in the sequel,
clearly $\hat V^n$ is a continuous function on $[0,\infty)\!
\times\! [0,\infty)$ and $\hat V^n$ restricted to $C_n$ and $D_n$
belongs to $C^2(\bar C_n)$ and $C^2(\bar D_n)$ respectively.
Finally, since $b_n$ is convex, we know that $b_n(\varPhi^1)$ is a
continuous semimartingale. This shows that the change-of-variable
formula with local time on surfaces \cite[Theorem 2.1]{Pe-3} is
applicable to $\hat V^n$ composed with $\varPhi =
(\varPhi^1,\varPhi^2)$ and using \eqref{7.1} this gives
\begin{align} \h{0pc} \label{8.6}
e^{-\lambda t}\:\! \hat V^n(\varPhi_t) &= \hat V^n(\varPhi_0) + \int
_0^t e^{-\lambda s} \big( \LL_\varPhi \hat V^n \m \lambda \hat V^n
\big)(\varPhi_s)\, ds + \int_0^t e^{-\lambda s}\:\! \hat V^n_{\varphi
_1}(\varPhi_s)\:\! \mu\;\! \varPhi_s^1\, dB_s^1\\ \notag &\h{13pt}+
\int_0^t e^{-\lambda s}\:\! \hat V^n_{\varphi_2}(\varPhi_s)\:\! \mu
\;\! \varPhi_s^2\, dB_s^2 - \int_0^t e^{-\lambda s}\:\! \hat V^n
_{\varphi_2}(\varPhi_s-)\, d\ell_s^{b_n}(\varPhi) \\ \notag &=
\hat V^n(\varPhi_0) - \int_0^t e^{-\lambda s}\:\! L(\varPhi_s)
\;\! I(\varPhi_s\! \in\! C_n)\, ds + M_t^n - \int_0^t e^{-\lambda
s}\:\! \hat V_{\varphi_2}(\varPhi_s)\, d\ell_s^{b_n}(\varPhi)
\end{align}
where $M_t^n = \int_0^t e^{-\lambda s}\:\! \hat
V_{\varphi_1}(\varPhi_s)\:\! \mu\;\! \varPhi_s^1\;\! I(\varPhi_s\!
\in\! C_n)\, dB_s^1 + \int_0^t e^{-\lambda s}\:\! \hat
V_{\varphi_2}(\varPhi_s)\:\! \mu\;\! \varPhi_s^2\;\! I(\varPhi_s\!
\in\! C_n)\, dB_s^2$ is a conti- nuous martingale for $t \ge 0$ and
$\ell^{b_n}(\varPhi)$ is the local time of $\varPhi$ on the curve
$b_n$ given by
\begin{equation} \h{1pc} \label{8.7}
\ell_t^{b_n}(\varPhi) = \PP \text{-}\! \lim_{\eps \downarrow 0} \frac
{1}{2 \eps} \int_0^t\! I( -\eps < \varPhi_s^2 \m b_n(\varPhi_s^1) <
\eps ) \, d \langle \varPhi^2 \m b_n(\varPhi^1),\varPhi^2 \m b_n(
\varPhi^1) \rangle_s
\end{equation}
for $t \ge 0$. To gain control over the final term in \eqref{8.6},
note that the It\^o-Tanaka formula yields
\begin{align} \h{0pc} \label{8.8}
\big( b_n(\varPhi_t^1) \m \varPhi_t^2 \big)^+ &= \big( b_n(\varPhi_0^1)
\m \varPhi_0^2 \big)^+\! +\! \int_0^t I(b_n(\varPhi_s^1) \m \varPhi_s^2\!
>\! 0)\, d( b_n(\varPhi^1) \m \varPhi^2 )_s + \frac{1}{2}\;\! \ell_t^{b_n}
(\varPhi) \\ \notag &= \big( b_n(\varPhi_0^1) \m \varPhi_0^2 \big)^+\!
+\! \int_0^t I(b_n(\varPhi_s^1) \m \varPhi_s^2\! >\! 0)\;\! \big( b_n'
(\varPhi _s^1)\;\! d \varPhi_s^1 - d \varPhi_s^2 \big) \\ \notag
&\h{13pt}+ \frac{1}{2} \int_0^t I(b_n(\varPhi_s^1) \m \varPhi_s^2\!
>\! 0)\! \int_0^\infty d \ell_s^{\psi_1}(\varPhi^1)\, d b_n'(\psi_1)
+ \frac{1}{2}\;\! \ell_t^{b_n}(\varPhi)
\end{align}
for $t \ge 0$ where $b_n'$ denotes the first derivative of $b_n$
whose existence follows by the implicit function theorem since
smooth fit fails at $b_n$ due to its suboptimality in the problem
\eqref{3.11}. Since $b_n$ is convex we see that $d b_n'$ defines a
non-negative measure on $[0,\infty)$ so that the double integral in
\eqref{8.8} is non-negative. It follows therefore from \eqref{8.8}
using \eqref{3.9}+\eqref{3.10} above that
\begin{align} \h{1pc} \label{8.9}
\frac{1}{2}\;\! \ell_t^{b_n}(\varPhi) &\le \big( b_n(\varPhi_t^1) \m
\varPhi_t^2 \big)^+ \m \int_0^t I(b_n(\varPhi_s^1) \m \varPhi_s^2\!
>\! 0)\;\! b_n' (\varPhi _s^1)\;\! \lambda (1 \p \varPhi_s^1)\, ds
\\ \notag &\h{14pt} +\! \int_0^t I(b_n(\varPhi_s^1) \m \varPhi_s^2\!
>\! 0)\;\! \lambda (1 \p \varPhi_s^2)\, ds + N_t^n
\end{align}
where $N_t^n = -\int_0^t I(b_n(\varPhi_s^1) \m \varPhi_s^2\! >\!
0)\;\! b_n' (\varPhi _s^1)\;\! \mu\;\! \varPhi_s^1\, dB_s^1 +
\int_0^t I(b_n(\varPhi_s^1) \m \varPhi_s^2\! >\! 0)\;\! \mu\;\!
\varPhi_s^2\, dB_s^2$ is a conti- nuous local martingale for $t \ge
0$. Let $(\tau_m)_{m \ge 1}$ be a localisation sequence of stopping
times for $N^n$, define the stopping time
\begin{equation} \h{7pc} \label{8.10}
\sigma_m = \inf\, \{\, t \ge 0\; \vert\; \varPhi_t^1 \le \tfrac{1}{m}
\, \}
\end{equation}
and set $\rho_m := \tau_m \wedge \sigma_m$ for $m \ge 1$. From
\eqref{8.9} we then find that
\begin{align} \h{3pc} \label{8.11}
\frac{1}{2}\;\! \EE_{\varphi_1,\varphi_2}^\infty \big[ \ell_{t\wedge
\rho_m}^{b_n} (\varPhi) \big] &\le \varphi_2^* - b_n' \big( \tfrac{1}
{m} \big)\! \int_0^t \lambda \big( 1 \p \EE_{\varphi_1,\varphi_2}^\infty
(\varPhi_s^1) \big)\, ds \\ \notag &\h{13pt}+ \int_0^t \lambda \big(
1 \p \EE_{\varphi_1,\varphi_2}^\infty (\varPhi_s^2) \big)\, ds \le
K_m(t)
\end{align}
for $t \ge 0$ and $m \ge 1$ where the positive constant $K_m(t)$
does not depend on $n \ge 1$ because each $b_n$ is convex and $b_n
\uparrow b$ on $[0,1/m]$ as $n \rightarrow \infty$ so that
$b_n'(1/m)$ must remain bounded from below over $n \ge 1$ if $b_n$
is to stay below $b$ on $[0,1/m]$ for all $n \ge 1$. In addition, by
\eqref{7.7} we know that $\hat V_{\varphi_2}$ is continuous on $\bar
C$ and hence uniformly continuous too because $\bar C$ is a compact
set. It follows therefore that $0 \le \hat
V_{\varphi_2}(\varphi_1,b_n(\varphi_1)) \le \eps$ for all $\varphi_1
\in [0,\varphi_0^n]$ and all $n \ge n_\eps$ with $n_\eps \ge 1$
large enough depending on the given and fixed $\eps>0$. Combining
this fact with \eqref{8.11}, upon replacing $t$ with $t \wedge
\rho_m$ in the final integral of \eqref{8.6} and taking
$\EE_\varphi^\infty$ of the resulting expression for $\varphi \in
[0,\infty)\! \times\! [0,\infty)$ given and fixed, we see that
\begin{equation} \h{3pc} \label{8.12}
0 \le \EE_\varphi^\infty \bigg[ \int_0^{t \wedge \rho_m}\!\!
e^{-\lambda s}\:\! \hat V_{\varphi_2}(\varPhi_s)\, d\ell_s^
{b_n}(\varPhi) \bigg] \le 2\:\! \eps\:\! K_t(m)
\end{equation}
for all $n \ge n_\eps$ with $t \ge 0$ and $m \ge 1$ given and fixed.
This shows that the expectation in \eqref{8.12} tends to zero as $n$
tends to infinity for every $t \ge 0$ and $m \ge 1$ given and fixed.
Using this fact in \eqref{8.6} upon replacing $t$ with $t \wedge
\rho_m$, taking $\EE_\varphi^\infty$ on both sides, and letting $n$
tend to infinity, we find by the monotone convergence theorem upon
recalling \eqref{7.6} that
\begin{equation} \h{3pc} \label{8.13}
\hat V(\varphi) = \EE_\varphi^\infty \big[ e^{-\lambda(t \wedge \rho
_m)} \:\! \hat V(\varPhi_{t \wedge \rho_m}) \big] + \EE_\varphi^\infty
\bigg[\! \int_0^{t \wedge \rho_m}\!\! e^{-\lambda s}\:\! L(\varPhi_s)
\:\! I(\varPhi_s\! \in\! C) \, ds \bigg]
\end{equation}
for all $t \ge 0$ and all $m \ge 1$. Letting $m \rightarrow \infty$
and using that $\rho_m \rightarrow \infty$ because $0$ is a natural
boundary point for $\varPhi^1$, we see from \eqref{8.13} upon
recalling \eqref{7.6} and using the dominated convergence theorem
that
\begin{equation} \h{3pc} \label{8.14}
\hat V(\varphi) = \EE_\varphi^\infty \big[ e^{-\lambda t}\:\! \hat
V(\varPhi_t) \big] + \EE_\varphi^\infty \bigg[\! \int_0^t e^{-\lambda
s} \:\! L(\varPhi _s)\:\! I(\varPhi_s\! \in\! C) \, ds \bigg]
\end{equation}
for all $t \ge 0$. Finally, letting $t \rightarrow \infty$ in
\eqref{8.14} and using the dominated and monotone convergence
theorems upon recalling \eqref{7.6}, we find that
\begin{equation} \h{5pc} \label{8.15}
\hat V(\varphi) = \EE_\varphi^\infty \bigg[\! \int_0^\infty e^{-\lambda
s}\:\! L(\varPhi _s)\:\! I(\varPhi_s\! \in\! C) \, ds \bigg]
\end{equation}
for all $\varphi \in [0,\infty)\! \times\! [0,\infty)$. Recalling
\eqref{6.10} and \eqref{8.2} above we see that this establishes the
representation \eqref{8.4} as claimed. Moreover, the fact that
$\tau_b$ from \eqref{8.5} is optimal in \eqref{3.11} follows by
\eqref{6.10} above. Finally, inserting $\varphi_2 = b(\varphi_1)$ in
\eqref{8.4} and using that $\hat V(\varphi_1,b(\varphi_1))=0$, we
see that $b$ solves \eqref{8.3} as claimed.

\vspace{4mm}

2.\ \emph{Uniqueness}. To show that $b$ is a unique solution to the
equation \eqref{8.3} in the specified class of functions, one can
adopt the four-step procedure from the proof of uniqueness given in
\cite[Theorem 4.1]{DuPe} extending and further refining the original
uniqueness arguments from \cite[Theorem 3.1]{Pe-1}. Given that the
present setting creates no additional difficulties we will omit
further details of this verification and this completes the proof.
\hfill $\square$

\vspace{4mm}

The nonlinear Fredholm integral equation \eqref{8.3} can be used to
find the optimal stopping boundary $b$ numerically (using Picard
iteration). Inserting this $b$ into \eqref{8.4} we also obtain a
closed form expression for the value function $\hat V$. Collecting
the results derived throughout the paper we now disclose the
solution to the initial problem.

\vspace{4mm}

\textbf{Corollary 6.} \emph{The value function in the initial problem
\eqref{2.5} is given by
\begin{equation} \h{5pc} \label{8.16}
V(\pi) = (1 \m \pi)\, \Big[\:\! 1 + c\;\! \hat V \Big( \frac{\pi}{1
\m \pi}, \frac{\pi}{1 \m \pi} \Big)\:\! \Big]
\end{equation}
for $\pi \in [0,1]$ where the function $\hat V$ is given by
\eqref{8.4} above. The optimal stopping time in the initial problem
\eqref{2.5} is given by
\begin{align} \h{2pc} \label{8.17}
\tau_* = \inf\, \Big \{\, t \ge 0\; &\big \vert\; e^{\;\! \mu X_t^2
+ (\lambda - \frac{\mu^2}{2})\:\! t} \Big( \frac{\pi}{1 \m \pi} +
\lambda \! \int_0^t\! e^{\;\! -\mu X_s^2 - (\lambda - \frac{\mu^2}
{2})\:\! s}\, ds \Big) \\ \notag &\h{6pt}\ge b \Big( e^{\;\! \mu X_t^1
+ (\lambda - \frac{\mu^2}{2})\:\! t} \Big( \frac{\pi}{1 \m \pi} +
\lambda \! \int_0^t\! e^{\;\! -\mu X_s^1 - (\lambda - \frac{\mu^2}
{2})\:\! s}\, ds \Big) \Big)\, \Big \}
\end{align}
where $b$ is a unique solution to \eqref{8.3} above (see Figure 1).}

\vspace{6mm}

\textbf{Proof.} The identity \eqref{8.16} was established in
\eqref{3.1} above. The explicit form of the optimal stopping time
\eqref{8.17} follows from \eqref{8.5} in Theorem 5 combined with
\eqref{2.14}+\eqref{2.15} above. The final claim on $b$ was derived
in Theorem 5 above. This completes the proof. \hfill $\square$

\section{Higher dimensions}

The quickest detection problem formulated in Section 2 and the
results derived in Sections 3-4 and 6-8 extend in a straightforward
way from dimension \emph{two} to dimension \emph{three or higher}.
This is readily obtained by replacing the coordinate number two of
the observed process $X$ by the coordinate number three or higher
throughout and only the notation gets more complicated. In this
section we briefly highlight this extension for future reference.

In the more general case, we consider a Bayesian formulation of the
problem \eqref{2.5} where it is assumed that one observes a sample
path of the standard $n$\!-dimensional Brownian motion $X=(X^1,
\ldots, X^n)$, whose coordinate processes $X^1, \ldots, X^n$ are
standard Brownian motions with zero drift initially, and then at
some random/unobservable time $\theta$ taking value $0$ with
probability $\pi \in [0,1]$ and being exponentially distributed with
parameter $\lambda>0$ given that $\theta>0$, one of the coordinate
processes $X^1, \ldots, X^n$ gets a (known) non-zero drift $\mu$
permanently. The problem is to detect the time $\theta$ at which a
coordinate process gets the drift $\mu$ as accurately as possible
(neither too early nor too late).

\vspace{6mm}

\textbf{Remark 7 (Higher dimensions).} All the results and arguments
in Sections 2-4 and 6-8 extend in an obvious way and remain valid
when the coordinate number $n$ is three or higher. The optimal
stopping boundary $b$ is no longer a curve but a surface in
$[0,\infty)^n$ which is obtained by replacing $b(\varphi_1)$ by
$b(\varphi_1, \ldots, \varphi_{n-1})$ above. In particular, the
existence and uniqueness results of Theorems 4 and 5 remain valid
when $n$ is three or higher and so does the solution to the initial
problem \eqref{2.5} as discussed in Corollary 6 above.

\vspace{6mm}

\textbf{Remark 8 (Signal-to-noise ratio).} An interesting question
is what we gain, if anything, by observing \emph{all} coordinate
processes of $X = (X^1, \ldots, X^n)$ simultaneously in real time
instead of a particular/individual coordinate process only when $n
\ge 2$. It appears to be evident that observing a single coordinate
process of one's choice is suboptimal, if for nothing else, then
because the drift $\mu$ may not appear in the chosen coordinate
process at all. Moreover, even if this deficiency is removed by
adding all coordinate processes and forming $Y_t := X_t^1 \p \ldots
\p X_t^n$ as the observed one-dimensional process for $t \ge 0$, we
see from \eqref{2.1} that $Y$ solves
\begin{equation} \h{7pc} \label{8.18}
dY_t = \mu\:\! I(t\! \ge\! \theta)\, dt + \sqrt{n}\;\! dW_t
\end{equation}
where $W_t := (\sum_{i=1}^n B_t^i)/\sqrt{n}$ is a standard Brownian
motion for $t \ge 0$. Comparing \eqref{8.18} with either \eqref{2.1}
or \eqref{2.2} we see that the \emph{signal-to-noise ratio} (defined
as the difference between the new drift and the old drift divided by
the diffusion coefficient) has \emph{decreased} in \eqref{8.18}
because $\mu/\sqrt{n} < \mu$ when $\mu$ is positive and $n \ge 2$.
Similarly, setting $Z_t := (\sum_{i=1}^n X_t^i)/\sqrt{n}$ for $t \ge
0$ we see from \eqref{8.18} that $Z$ solves
\begin{equation} \h{7pc} \label{8.19}
dZ_t = \frac{\mu}{\sqrt{n}}\:\! I(t\! \ge\! \theta)\, dt + dW_t\, .
\end{equation}
Thus, assuming that $Z_t$ is being observed for $t \ge 0$, we see
that the quickest detection problem for $Z$ reduces to the problem
in one dimension considered in Section 5 above. From \eqref{8.19} we
see however that the drift $\mu/\sqrt{n}$ in the former problem is
\emph{strictly smaller} that the drift $\mu$ in the latter problem
when $n \ge 2$ so that quickest detection for the observed process
$Z$ is \emph{harder}. The final result of Corollary 6 above
(combined with Remark 7) shows that quickest detection of a
coordinate drift requires full knowledge of all coordinate processes
$X^1, \ldots, X^n$, so that observing one of them only, or even
their sum, is insufficient to reach full optimality.

\vspace{4mm}

\textbf{Acknowledgements.} The authors gratefully acknowledge
support from the United States Army Research Office Grant
ARO-YIP-71636-MA.

\begin{center}

\end{center}


\par \leftskip=24pt

\vspace{10mm}

\ni Philip A.\ Ernst \\
Department of Statistics \\
Rice University \\
6100 Main Street \\
Houston TX 77005 \\
United States \\
\texttt{philip.ernst@rice.edu}

\leftskip=25pc \vspace{-36mm}

\ni Goran Peskir \\
Department of Mathematics \\
The University of Manchester \\
Oxford Road \\
Manchester M13 9PL \\
United Kingdom \\
\texttt{goran@maths.man.ac.uk}

\par

\end{document}